\documentclass[pra, preprint, superscriptaddress]{revtex4}
\usepackage{graphicx}
\usepackage{mathrsfs}
\usepackage{dcolumn}
\usepackage{amsmath}
\usepackage{amssymb}
\usepackage{color}
\usepackage{float}
\usepackage{epstopdf}

\def\beq{\begin{equation}}
\def\eeq{\end{equation}}
\def\eq#1{(\ref{#1})}
\def\simg{\,\hbox{\kern.1em \lower.6ex \hbox{$\sim$} \kern-1.12em
          \raise.6ex \hbox{$>$} }}
\def\siml{\,\hbox{\kern.1em \lower.6ex \hbox{$\sim$} \kern-1.12em
          \raise.6ex \hbox{$<$} }}

\def\XXint#1#2#3{{\setbox0=\hbox{$#1{#2#3}{\int}$}
     \vcenter{\hbox{$#2#3$}}\kern-.5\wd0}}

\newcommand\bea{\begin{eqnarray}}
\newcommand\eea{\end{eqnarray}}

\begin{document}

\title{On the asymptotic distinct prime partitions of integers}
\author{M. V. N. Murthy}
\affiliation{The Institute of Mathematical Sciences, Chennai 600 113, India}
\email{murthy@imsc.res.in}
\author{Matthias Brack}
\affiliation{Institute of Theoretical Physics, University of Regensburg,
D-93040 Regensburg, Germany}
\email{matthias.brack@ur.de}
\author{R. K. Bhaduri\footnote{Deceased in November 2019 after the 
submission of this paper}}
\affiliation{Department of Physics and Astronomy, McMaster University,
Hamilton L8S4M1, Canada}

\begin{abstract}

We discuss $Q(n)$, the number of ways a given integer $n$ may be written 
as a sum of distinct primes, and study its asymptotic form $Q_{as}(n)$ 
valid in the limit $n\to\infty$. We obtain $Q_{as}(n)$ by Laplace 
inverting the fermionic partition function of primes, in number theory called 
the generating function of the distinct prime partitions, in the saddle-point 
approximation. We find that our result of $Q_{as}(n)$, which includes two
higher-order corrections to the leading term in its exponent and a pre-exponential
correction factor, approximates the exact $Q(n)$ far better than 
its simple leading-order exponential form given so far in the literature. 

\bigskip

\end{abstract}

\maketitle

\section{Introduction}
\label{secint}

The asymptotic form $Q_{as}(n)$ of the distinct prime partitions $Q(n)$, 
called the sequence A000586 in \cite{oeis}, has in the literature so far 
only been given in its leading exponential form \cite{roth}, published 
as long ago as in 1954. We believe that it is time to improve the 
asymptotics by corrections to the leading exponential form which itself 
is rather poor \cite{oeis}. Although this might appear as a 
straight-forward task, it requires quite some cumbersome algebraic 
efforts. For the unrestricted prime partions $P(n)$, called the sequence 
A000607 in \cite{oeis}, the next-to-leading order correction has been 
derived and published in 2008 in this journal \cite{vaughan}. The 
numerical coefficient of the next-to-leading order term was corrected 
recently in \cite{BBBM}. For our present investigations we will employ 
here the same method as in \cite{BBBM}. Although this method 
algebraically is very close to that used in \cite{vaughan}, we present 
it as viewed from the standpoint of statistical mechanics. The fact 
that the two partitions mentioned above can be connected to 
many-body systems of particles obeying opposite statistics 
may, in our opinion, serve as a bridge between mathematical physics and 
pure number theory.

It is well established by now \cite{muoi} that the techniques of 
statistical mechanics can be applied to obtain any type of partition of 
a positive integer $n$. The partition function of a gas in statistical 
mechanics contains information on the distribution of the total energy 
among the constituents and hence plays the same role as the generating 
function of the corresponding partitions in number theory. This relation 
was used in \cite{muoi}, where the number of partitions $P(n)$ known 
from number theory is obtained from the quantum density of states 
$\rho(E)$ given by the inverse Laplace transform of the partition 
function. Taking the Laplace transform approximately using the 
saddle-point method then yields the asymptotic forms $P_{as}(n)$.

This method was applied in \cite{BBBM} to the unrestricted partitions 
$P(n)$ of integers into primes, i.e., the sequence A000607. For a system 
whose single-particle levels are defined by the primes $p$ as an ordered 
set, the total energy is given by a sum of primes, and the corresponding 
density of states is related to the {\it number of unrestricted prime 
partitions} $P(n)$, assuming that the particles behave like {\it 
bosons}. The asymptotic form $P_{as}(n)$ obtained in \cite{BBBM} was 
found to approximate the exact $P(n)$ for large $n$ much better than the 
asymptotic expressions given earlier in the literature 
\cite{yang,vaughan}. The same method was also applied more recently to 
distinct square partitions in \cite{MBBB}.

In the present paper we study $Q(n)$, the number of ways a given integer $n$ may 
be written as a sum of {\it distinct primes}, i.e., the sequence A000586. As 
only distinct primes are allowed in $Q(n)$, this corresponds to a system of 
{\it fermionic} particles, obeying the Pauli exclusion principle, still with 
the primes $p$ as single-particle levels. We again use the saddle-point method 
for Laplace inverting their partition function to derive algebraically the 
asymptotic form $Q_{as}(n)$ valid in the limit $n\to\infty$. In numerical 
computations up to $n=10^5$ we find that, like in \cite{BBBM}, our result for 
$Q_{as}(n)$ approaches the exact $Q(n)$ far better than the simple leading-order 
exponential form given so far in the literature \cite{roth}.

The plan of our paper is as follows. In Section \ref{secpf}, we establish the 
relation of $Q(n)$ to the partition function and the density of states 
$\rho^F(E)$. In Sec.\ \ref{secasypf} we derive its asymptotic form using the 
saddle-point method, and in Sec.\ \ref{secspeq} we give the explicit solution of 
the saddle-point equation leading to our final result for $Q_{as}(n)$. 
In Sec.\ \ref{secnum} our asymptotic result is compared numerically with the 
exact function $Q(n)$ for the distinct prime partitions. We conclude the 
paper with a short summary in Sec.\ \ref{secsum}.

\section{Partitions into primes}

\subsection{Fermionic partition function and its relation to $Q(n)$}
\label{secpf}

Consider a large number $N$ of fermions whose single-particle spectrum
is given by the primes $p$. The total energy $E$ of the system is given by 
\beq 
E=\sum_{p} n_p\, p\, .
\label{Etot}
\eeq 
(We use throughout dimensionless variables and take the particle mass $m$, 
the Planck constant $\hbar$ and the Boltzmann constant $k$ to be unity: 
$m=\hbar=k=1$.)
Here and in the following, the sums $\sum_p$ run over all primes $p$, and $n_p$ 
are the fermionic occupancies of the levels which must be zero or one, such that
\beq
\sum_p n_p = N\,, \qquad n_p=0,1.
\label{Nsum}
\eeq 
The number of possible energy partitions $E$ with the restriction \eq{Nsum} 
shall be denoted by $Q_N(E)$, where the subscript $N$ keeps track of 
the total number of particles. Although $E$ is necessarily integer, we treat it 
as a continuous variable like in statistical mechanics. $Q_N(n)$ is the
number of $N$-{\it restricted} fermionic partitions of $n$, i.e., the number of
ways to write $n$ as a {\it sum of $N$ distinct primes.} In the limit 
$N\to\infty$, $Q_N(n)$ will tend towards the number of {\it unrestricted} 
but {\it distinct prime partitions} $Q(n)$ under consideration here.

For the purpose of this paper, we are only interested in the limit $N\to\infty$ 
of the fermionic partitions which then become unrestricted as stated above. The 
{\it quantum-statistical partition function} $Z^F(\beta)$ is in this limit given by
\beq
Z^F(\beta)=\prod_p [1+e^{-\beta\,p}]
\label{zinfty}
\eeq
where $\beta=1/kT$ is the inverse temperature and the product runs over all 
primes $p$. Taylor expanding the expontential in \eq{zinfty} and reordering the 
terms yields the alternative form of the partition function
\beq
Z^F(\beta)=\sum_{n=0}^\infty Q(n)\,e^{-n\beta}\,,
\label{zgen}
\eeq
which in number theory is known as the {\it generating function} of the $Q(n)$.
In the On-line Encyclopedia of Integer Sequences (OEIS) \cite{oeis}, the
sequence of numbers $Q(n)$ is called the sequence A000586. Its first
ten members are $Q(n)$ = 1, 0, 1, 1, 0, 2, 0, 2, 1, 1 for $n=0,...,9$,
where $Q(0)=1$ by definition. Note also that the $Q(n)$ are a subset of the 
(bosonic) prime partitions $P(n)$, called the sequence A000607 in \cite{oeis}.

From the partition function, we obtain the {\it many-body density of states} 
$\rho^F(E)$ by an inverse Laplace transform:
\beq
\rho^F(E) = \frac{1}{2\pi i}\int_{-i\infty}^{i\infty} d\beta\, 
              Z^F(\beta)\,e^{\beta E}.
\label{rhodef}
\eeq
Hereby $\beta$ is taken as a complex variable and the integration above runs
along the imaginary axis of the complex $\beta$ plane. Later the Laplace 
inversion shall be taken in the saddle-point approximation.

It is important now to realize that $Q(n)$ is related to the density of
states $\rho^F(E)$ in the following way. Taking directly the exact inverse Laplace 
transform of \eq{zgen}, we find
\beq
\rho^F(E)=\sum_{n=0}^\infty Q(n)\,\delta(E-n)\,,
\label{denq}
\eeq
where $\delta(E-n)$ is the Dirac delta function peaked at $E=n$. We see thus 
that $\rho^F(E)$ can also be understood as the {\it density of distinct prime 
partitions}. Like it was argued in \cite{MBBB} for the distinct square
partitions, averaging $\rho^F(E)$ over a sufficiently large energy interval
$\Delta E$ is asymptotically the same as averaging $Q(n)$ over a 
sufficiently large interval $\Delta n$:
\beq
\langle\, \rho^F(E)\, \rangle_{\Delta E} \; 
          \sim \; \langle\, Q(n)\, \rangle_{\Delta n}
          \quad \hbox{ for } \; E,n \rightarrow \infty\,.
\label{avs}
\eeq 
Therefore determing the {\it asymptotic average part} $\rho_{as}^F(E)$ of the density 
of states valid in the limit $E\to\infty$, which can be obtained by the 
saddle-point approximation to its inverse Laplace transform \eq{rhodef}, and 
equating $E=n$ will give the {\it average asymptotic form} $Q_{as}(n)$ of 
the distinct prime partitions.


\subsection{Asymptotic partition function from saddle-point approximation}
\label{secasypf}

We first rewrite the inverse Laplace transform \eq{rhodef} by taking the natural
log of $Z^F$ into the exponent:
\beq
\rho^F(E) = \frac{1}{2\pi i}\int_{-i\infty}^{+i\infty} d\beta\,
            e^{\beta E + \ln Z^F(\beta)}\,.
\label{invlap}
\eeq
We now evaluate this integral using the saddle-point method (also called the method 
of steepest descent). We define the exponent above as the canonical entropy function
\beq
S^F(E,\beta)=\beta E + \ln Z^F(\beta)\,.
\label{sdef}
\eeq
Applying the saddle-point method to \eq{invlap} requires to find a stationary point 
$\beta_0$ of the function $S^F(E,\beta)$ by solving the {\it saddle-point equation} 
\beq
\left. \frac{\partial S^F(E,\beta)}{\partial \beta}\right|_{\beta_0}
     = E+\frac{\partial Z^{F}(E,\beta_0)/\partial\beta}{Z^F(E,\beta_0)}=0\,.
\label{spc}
\eeq
If this equation has a solution $\beta_0$, which will be a function $\beta_0(E)$, 
one evaluates the successive partial derivatives of $S^F(E,\beta)$ at $\beta_0$:
\beq
S^{F(n)}(E,\beta_0) = \left. \frac{\partial^n S^F(E,\beta)}{\partial \beta^n}
                    \right|_{\beta_0}.
\eeq 
The approximate result of the inverse Laplace transform then is given by
\beq
\rho_{as}^F(E) = \frac{e^{S^F(E,\beta_0)}}{\sqrt{2\pi S^{F(2)}(E,\beta_0)}}
                 \left[1+\cdots\right],
\label{rhosol}
\eeq
where the dots indicate the so-called cumulants involving higher derivatives of the 
entropy, which become more important for large $\beta$ (see, e.g., Ref.\ 
\cite{jelovic}). Since we are interested here in the limit $\beta\to 0$ relevant 
for the asymptotics of large $E$, we can neglect these cumulants. 

Next, we take the natural log of the partition function given in Eq.(\ref{zinfty})
\begin{equation}
\ln Z^F(\beta) = \sum_{p=2}^{\infty}\ln\left(1+e^{-\beta p}\right). 
\label{lnzbex}
\end{equation}
and approximate it by the integral 
\begin{equation}  
\ln Z^F(\beta) \sim \int_{2}^{\infty}dx\, g_{av}(x) \ln\left(1+e^{-\beta x }\right), 
\label{lnzb}
\end{equation}
where $g_{av}(x)=\frac{1}{\ln(x)}$ is the approximate density of primes,
using the prime number theorem. If the density $g_{av}(x)$ were exact, then 
the integral would give the exact result \eq{lnzbex}. 

The evaluation the integral in the limit $\beta\to 0$ follows closely the 
method outlined in \cite{BBBM}. Denoting $y=\beta x$, the integral becomes
\begin{equation}  
\ln Z^F(\beta) \sim \frac{1}{\beta}\int_{2\beta}^{\infty} dy\,
                    \frac{1}{\ln(\frac{y}{\beta})}\ln\left(1+e^{-y }\right)   
                  = -\frac{1}{\beta\ln\beta}\int_{2\beta}^{\infty} dy\,
                    \frac{1}{1-\frac{\ln(y)}{\ln(\beta)}}\ln\left(1+e^{-y }\right).
\end{equation}
In the limit $\beta\to 0$ we may write this integral
as an asymptotic series
\begin{equation}
\ln Z^F(\beta) \sim  -\frac{1}{\beta\ln\beta}\int_{2\beta}^{\infty} dy
                     \left[1+\sum_{k=1}^\infty 
                     \left(\frac{\ln(y)}{\ln(\beta)}\right)^{\!k}\right]
                     \ln\left(1+e^{-y }\right).
\end{equation}
This is now a series in the expansion parameter $1/\ln(\beta)$ since each
term is divided by the power $(\ln\beta)^k$. As we shall see later, in
the leading saddle-point approximation $\ln(\beta_0)\approx \ln(E)$ and hence
this is an asymptotic series in $1/\ln(E)$ as well. In the asymptotic limit
we take the lower limit of the integral to be zero. 
 
For the present analysis, we retain the leading term and the first correction, 
like for the bosonic prime partitions in \cite{BBBM}, and define
\begin{equation}  
\ln Z^F_{as}(\beta) = -\frac{1}{\beta\ln\beta}\int_{0}^{\infty} dy
                     \left[1+\left(\frac{\ln(y)}{\ln(\beta)}\right)\right]
                     \ln\left(1+e^{-y }\right).
\end{equation}
The integrals may again be evaluated analytically and we obtain
\begin{equation}  
\ln Z^F_{as}(\beta) = \frac{1}{\beta\ln(\beta)}\left[-\frac{\pi^2}{12}
                     +\frac{1}{\ln\beta}\left(\frac{C\pi^2}{12}
                     +\sum_k(-1)^{k-1}\frac{\ln(k)}{k^2}\right)\right]\!,
\label{lnzas}
\end{equation}
where $C=0.5772156649\!\dots$ is the Euler constant.

\subsection{Solution of saddle-point equation and  $Q_{as}(n)$}
\label{secspeq}

In order to find the saddle point $\beta_0$ from Eq.\ \eq{spc}, we start from 
the entropy $S^F(\beta)$ in the asymptotic limit. Using Eqs.\ (\ref{sdef}) and 
(\ref{lnzas}) we get up to order $1/(\ln\beta)^2$
\begin{equation}  
S^F(E,\beta) = \beta E-\frac{F_1}{\beta\ln(\beta)}
               +\frac{F_2}{\beta\,(\ln\beta)^2}+\cdots
\label{sfbeta}
\end{equation}
where 
\begin{equation}  
F_1 = \frac{\pi^2}{12}\,, \qquad 
F_2=\left[\frac{C\pi^2}{12}+\sum_{k=1}^\infty (-1)^{k-1}
                                 \frac{\ln(k)}{k^2}\right]
\end{equation}
The infinite sum in $F_2$ may be expressed in a closed form in terms of 
a derivative of the Riemann zeta function, leading to
\begin{equation}
F_2=\frac{\pi^2}{12}[C-\ln(2)]-\frac{\zeta'(2)}{2}=0.3734242774\!\dots
\end{equation}

Eq.\ (\ref{sfbeta}) is identical in form with that of the bosonic case given 
in \cite{BBBM}:
\begin{equation}  
S^B(E,\beta)= \beta E-\frac{f_1}{\beta\ln(\beta)}
            + \frac{f_2}{\beta(\ln\beta)^2}+\cdots\,,
\label{sbbeta}
\end{equation}
where 
\begin{equation}  
f_1=\frac{\pi^2}{6}\,, \qquad f_2=\left[\frac{C\pi^2}{6}+\sum_k
                              \frac{\ln(k)}{k^2}\right]
\end{equation}
The constant $f_2$ may also be expressed in a closed form by
\begin{equation}
f_2=\frac{C\pi^2}{6}-\zeta'(2)=
\frac{\pi^2}{6}[12\ln(\!A)-\ln(2\pi)]=1.887029965\dots,
\end{equation}
where $A= 1.282427129100$... is the Glaisher-Kinkelin constant (see A074962 in \cite{oeis}).

The only difference in going from bosonic to the fermionic case is that 
the coefficients $f_1$ and $f_2$ of \cite {BBBM} are
replaced here by the $F_1$ and $F_2$, respectively. Therefore we obtain our 
result simply by replacing the coefficients $f_i$ in the bosonic case by the 
$F_i$ in the present fermionic case and following the steps outlined in \cite{BBBM}. 

Thus we can directly give the result for the fermionic case as
\begin{equation}
\rho_{as}^F(E) = \frac{1}{{\sqrt{4E^{3/2}[6\ln(E)]^{1/2}}}}\,
                \exp\left\{ 2\pi\sqrt{\frac{E}{6\ln(E)}}
                \left[1-\frac{1}{2}\frac{\ln[\ln(E)]}{\ln(E)}
                +b^F\frac{1}{\ln(E)}\right]\right\}
\end{equation}
with the constant
\beq
b^F = \left[\frac{F_2}{F_1}+\ln(\pi/\sqrt{6})\right]=0.7028796287\!\dots
\eeq
The asymptotic $Q_{as}(n)$ is then obtained replacing $E$ by $n$ above, so that:
\begin{equation}
Q_{as}(n)=\frac{1}{{\sqrt{4n^{3/2}[6\ln(n)]^{1/2}}}}\,
                \exp\left\{ 2\pi\sqrt{\frac{n}{6\ln(n)}}
                \left[1-\frac{1}{2}\frac{\ln[\ln(n)]}{\ln(n)}
                +b^F\frac{1}{\ln(n)}\right]\right\}.
\label{qasn} 
\end{equation}
This is the main result of the present paper. The corresponding result for the
bosonic partitions in \cite{BBBM} was
\begin{equation}
P_{as}(n)=\frac{1}{{\sqrt{4n^{3/2}[3\ln(n)]^{1/2}}}}\,
                \exp\left\{ 2\pi\sqrt{\frac{n}{3\ln(n)}}
                \left[1-\frac{1}{2}\frac{\ln[\ln(n)]}{\ln(n)}
                +b^B\frac{1}{\ln(n)}\right]\right\},
\label{pasn} 
\end{equation}
with the constant
\beq
b^B = \left[\frac{f_2}{f_1}+\ln(\pi/\sqrt{3})\right]=0.7426003995\!\dots
\eeq
Note that the leading exponential terms and the denominators of the pre-exponential 
terms in \eq{qasn} and \eq{pasn} differ by a factor $1/\!\sqrt{2}$. Note that 
since the $Q(n)$ are a subset of the $P(n)$, their values must be 
smaller, which asymptotically is brought about by the extra factor $1/\!\sqrt{2}$ 
in the leading exponential term. The first correction term in the exponent, namely 
$-\frac{1}{2}\ln[\ln(n)]/\!\ln(n)$, is identical in both cases. As far as we know, 
the above result \eq{qasn} for the distinct prime partitions has not been given in 
the literature so far.

In the next section, we compare numerically our asymptotic result \eq{qasn} 
with the exact values $Q(n)$ of the distinct prime partitions.


\section{Numerical test of $Q_{as}(n)$}
\label{secnum}

\begin{figure}[h]
\centering
\vspace*{-0.3cm}
\includegraphics[width=12.cm,angle=0]{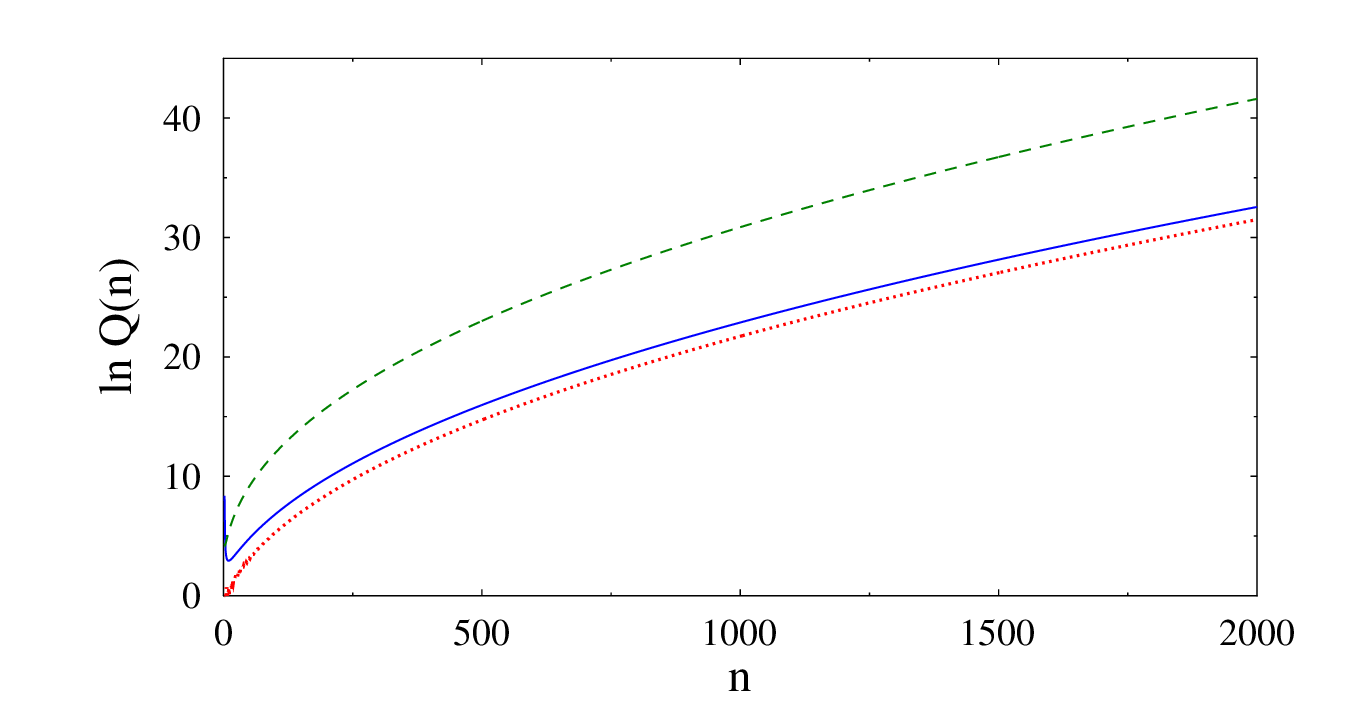}
\vspace*{-0.7cm}
\caption{Exact $\ln Q(n)$ by dots (red), lowest-order asymptotic form 
$\ln Q_0(n)$ from \eq{qas0} by the dashed line (green) and our full asymtotic 
form $\ln Q_{as}(n)$ given by \eq{qasn} by the solid (blue) line, shown as 
functions of $n$ up to $n=2000$.}
\label{figqofnlow}
\end{figure}

\begin{figure}[h]
\centering
\vspace*{-0.3cm}
\includegraphics[width=12.cm,angle=0]{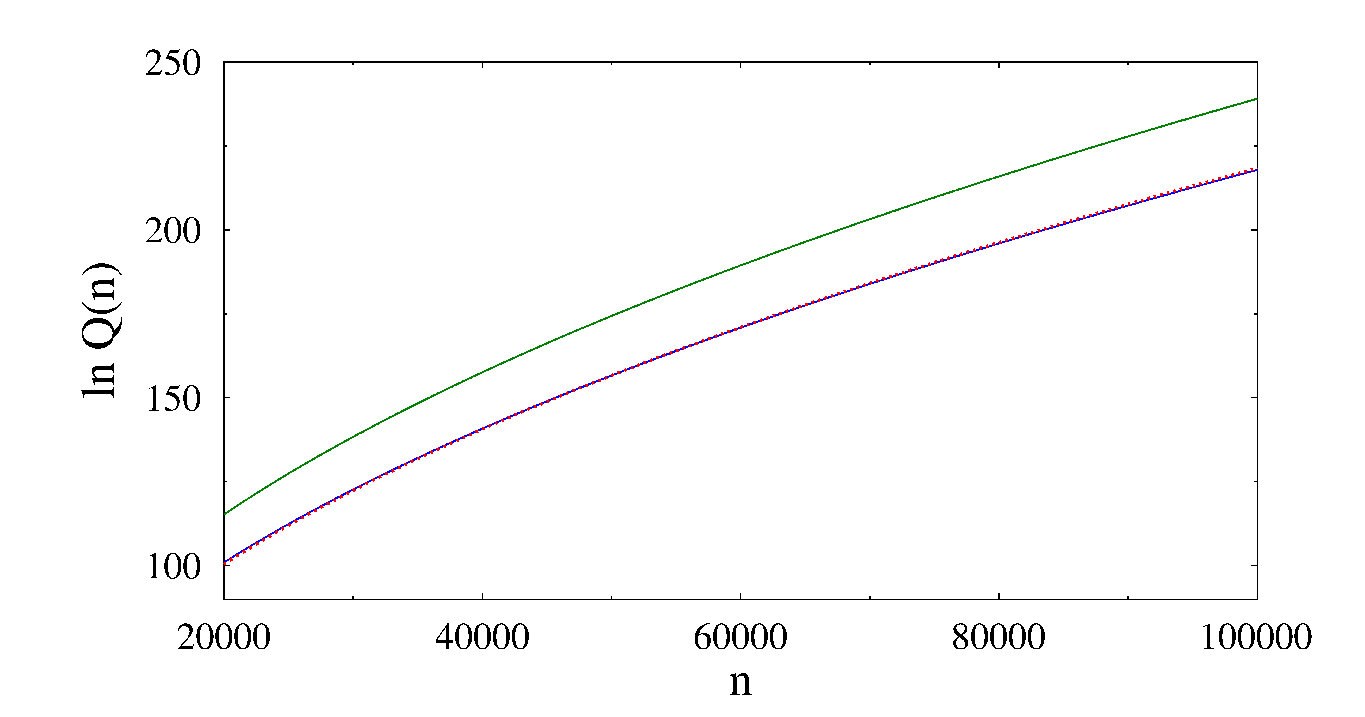}
\vspace*{-0.7cm}
\caption{Same as Fig.\ \ref{figqofnlow} in the region $n=20\,000-100\,000$.}
\label{figqofn}
\end{figure}

In this section we test our asymptotic result \eq{qasn} numerically. We have 
generated the exact $Q(n)$ up to $n=100\,000$. In Figs.\ \ref{figqofnlow}
and \ref{figqofn} we show their values by the dots (red) on a logarithmic scale 
in two regions of $n$. The dashed line (green) shows the leading-order exponential 
expression  
\beq
Q_0(n)=\exp\left\{2\pi\sqrt{\frac{n}{6\ln(n)}}\right\}\,,
\label{qas0} 
\eeq
which is the only asymptotic result that has been given so far in the literature
\cite{roth}, and the solid (blue) line gives our full asymptotic result \eq{qasn}. 

\noindent
A large discrepancy between $Q_0(n)$ and $Q(n)$ is noticed for all
$n$. Our full asymptotic result $Q_{as}(n)$ \eq{qasn} approaches the exact 
$Q(n)$ much better (except in the academic limit $n\to 0$ where it 
diverges due to the pre-exponential factor). In Fig.\ \ref{figqofn} for the
values $n \geq 20\,000$, the two curves can hardly be distinguished. 

We have thus achieved a considerable improvement over the simple exponential form 
\eq{qas0}. A closer look reveals that the curve for $Q_{as}(n)$, which for smaller 
$n$ overestimates the exact $Q(n)$, crosses the curve of the latter around
$n\sim 50\,000$. A similar result was found in \cite{BBBM} for the bosonic 
prime partitions, where $P_{as}(n)$ crosses  $P(n)$ much earlier and then tends
to approach it asymptotically from below for $n\to\infty$. 

\newpage

\begin{figure}[H]
\centering
\vspace*{-0.3cm}
\includegraphics[width=12.cm,angle=0]{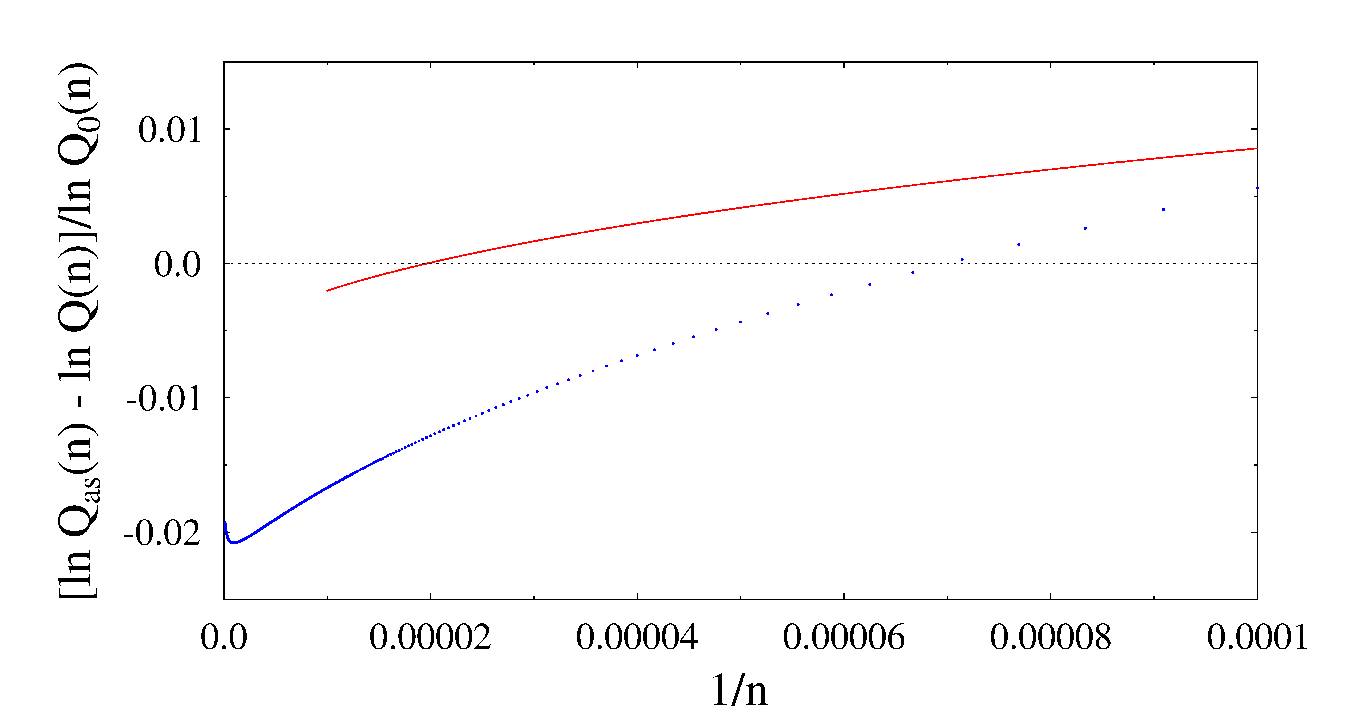}
\vspace*{-0.2cm}
\caption{Relative difference $[\ln Q_{as}(n)-\ln Q(n)]/\ln Q_0(n)$,
shown versus $1/n$ by the solid (red) line. The dotted (blue) line shows the
corresponding quantity obtained in \cite{BBBM} for the bosonic prime partitions 
$P(n)$.}
\label{figqperr}
\end{figure}

In order to focus on this asymptotic behavior, we show in Fig.\ \ref{figqperr} 
the difference of the natural logs relative to the lowest-order term, i.e., the 
quantity $[\ln Q_{as}(n)-\ln Q(n)]/\ln Q_0(n)$, plotted versus $1/n$ in a
region of the largest $n$ available. 
The solid (red) curve gives the result obtained with our full asymptotic form 
\eq{qasn}. For comparison we show in this figure by the dotted (blue) curve also 
the corresponding quantity obtained in Ref.\ \cite{BBBM} from the unrestricted
(bosonic) prime partitions $P(n)$ and their respective asymptotic forms. The 
overall behaviour of the two curves is similar. For the results in \cite{BBBM} 
we had larger values of $n$ available. There we noticed a tendency for the 
difference to approach zero from below for $1/n\to 0$ (i.e.\ $n\to\infty$), 
as can be recognized from the blue curve in Fig.\ \ref{figqperr}.

\noindent{\it Note added after completion of our work:}

\medskip

After the publication of our results on the arXiv server in 2019 \cite{MBBarx}, 
V. Kotesovec has performed numerical studies of our 
$\ln Q_{as}(n)$ and the exact $\ln Q(n)$, computing these quantities up to 
$n_{max}=10^8$. In order to achieve this, he programmed in the assembler a special 
floating point arithmetic in which both the mantissa and the exponent of these
quantities have 8 bytes.
By this procedure it was possible to generate 10 million terms in 9 hours, and the 
calculation of $n_{max}=10^8$ terms took 31 days. 

V.\ Kotesovec has kindly sent us his 
results which confirm our findings (red line) up to $n=10^5$ and furthermore show 
that, indeed, the difference $[\ln Q_{as}(n)-\ln Q(n)]/\!\ln Q_0(n)$ tends towards zero 
from the same side as that of the bosonic (unrestricted) prime partitions (blue curve). 
Similarly to the situation for the latter, the asymptotic ratio 
$\ln Q(n)/\!\ln Q_{as}(n)$ 
first exceeds the value 1 but then reaches a maximum, occurring here at 
$n = 14\,474\,250$, has an inflection point at $n \simg 33\,272\;000$, and gradually 
decreases back towards 1. One or two million terms are far from enough for this 
finding; it is necessary to have at least 40 million terms. A graph of Kotesovec's
result for $\ln Q(n)/\!\ln Q_{as}(n)$ is posted at OEIS \cite{kotes}.


\section{Summary}
\label{secsum}

In summary, we have shown how an improved asymptotic expression for the function 
$Q(n)$, which counts the number of distinct prime partitions of an integer $n$, 
can be obtained from asymptotic expansions of the partition function $Z^F(\beta)$ 
in \eq{zgen} and the corresponding density of states $\rho^F(E)$ in \eq{rhodef}. 
$Z^F(\beta)$ can be understood as the quantum-statistical partition function of a 
system of $N$ fermions, whose single-particle energy spectrum is given by the 
primes $p$, in the limit $N\to\infty$. It is identical to the generating function 
of the $Q(n)$ known in number theory. The density of states $\rho^F(E)$ is 
identical to the the density of distinct prime partitions given in Eq.\ \eq{denq}. 
Exploiting the connection between $\rho^F(E)$ and $Q(n)$ using the saddle-point 
approximation for the inverse Laplace transform \eq{rhodef}, we have obtained the
asymptotic form $Q_{as}(n)$ in Eq.\ \eq{qasn} and shown it numerically to approach 
the exact $Q(n)$ in the limit $n\to\infty$ far better than the hitherto known 
expression $Q_0(n)$ given in \eq{qas0}. 

We have used the same method as in Ref.\ \cite{BBBM} where the non-distinct 
prime partitions $P(n)$ were studied, and have found similar results as 
there. The asymptotic $Q_{as}(n)$ overestimates the exact $Q(n)$ for smaller 
$n$ but overshoots it for $n \simg 50,000$. Like in \cite{BBBM}, the 
limit $Q_{as}(n)\to Q(n)$ for $n\to\infty$ cannot be demonstrated rigorously.
However, forcing the calculation of our $Q_{as}(n)$ and of $Q(n)$ up to $n=10^8$, 
V.\ Kotesovec has shown numerically that $Q_{as}(n)$, indeed, approaches $Q(n)$ 
monotonously for $n \simg 4\times 10^7$ \cite{kotes}. He assumed that the difference
can be approximated by a term $c_2/\ln^2(n)$ in the square brackets of the 
asymptotic expansion \eq{qasn} and showed that the results are very sensitive
to the value of $c_2$. We join his suggestion that the systematic evaluation
of the algebraic value of $c_2$, or of other correction terms in \eq{qasn},
could be a topic of interesting future research for the next generation of
patient researchers.

But already now we can state that already with our present result \eq{qasn}, 
we have obtained an excellent asymptotic approximation for the distinct prime 
partitions which is far superior to the hitherto known result \cite{roth}.

\medskip

M.V.N.M. and M.B.\ acknowledge stimulating earlier correspondence with 
V. Kotesovec and, in particular, the communication of his most recent numerical
results. R.K.B.\ is grateful to the IMSc, Chennai, for its hospitality
during the final stages of our collaboration.


\newpage

\begin{thebibliography}{99}

\vspace*{-0.5cm}

\bibitem{oeis}     The On-line Encyclopedia of Integer Sequences (OEIS),
                   see $<$http://oeis.org$>$.

\bibitem{roth}     K. F. Roth and G. Szekeres, Q. J. Math. Oxf. Ser.\ 
                   (2) {\bf 5}, 241 (1954).
 
\bibitem{vaughan}  R. C. Vaughan, Ramanujan J. {\bf 15}, 109 (2008). 

\bibitem{BBBM}     J. Bartel, R. K. Bhaduri, M. Brack and M. V. N. Murthy, 
                   Phys.\ Rev.\ E {\bf 95}, 052108 (2017).

\bibitem{muoi}     M. N. Tran, M. V. N. Murthy and R. K. Bhaduri, 
                   Ann. Phys. {\bf 311}, 204 (2004).

\bibitem{yang}     Yifan Yang, Trans. Am. Math. Soc. {\bf 352}, 2581 (2000).

\bibitem{MBBB}     M. V. N. Murthy, M. Brack, R. K. Bhaduri, and J. Bartel, 
                   Phys. Rev.\ E {\bf 98}, 052131 (2018).

\bibitem{jelovic}  A. Jelovic, Phys. Rev. C {\bf 76}, 017301 (2007).

\bibitem{MBBarx}   M. V. N. Murthy, M. Brack, amd R. K. Bhaduri, 
                   $<$https://arxiv.org/abs/1904.02776$>$.
 
\bibitem{kotes}    Figure shown at $<$http://oeis.org/A000586$>$ via the link:
                   ``Vaclav Kotesovec, Plot log(a(n)) / log(Qas(n)) for 
                   $n = 2\,.\,.\, 10^8$''.













\end{thebibliography}
\end{document}